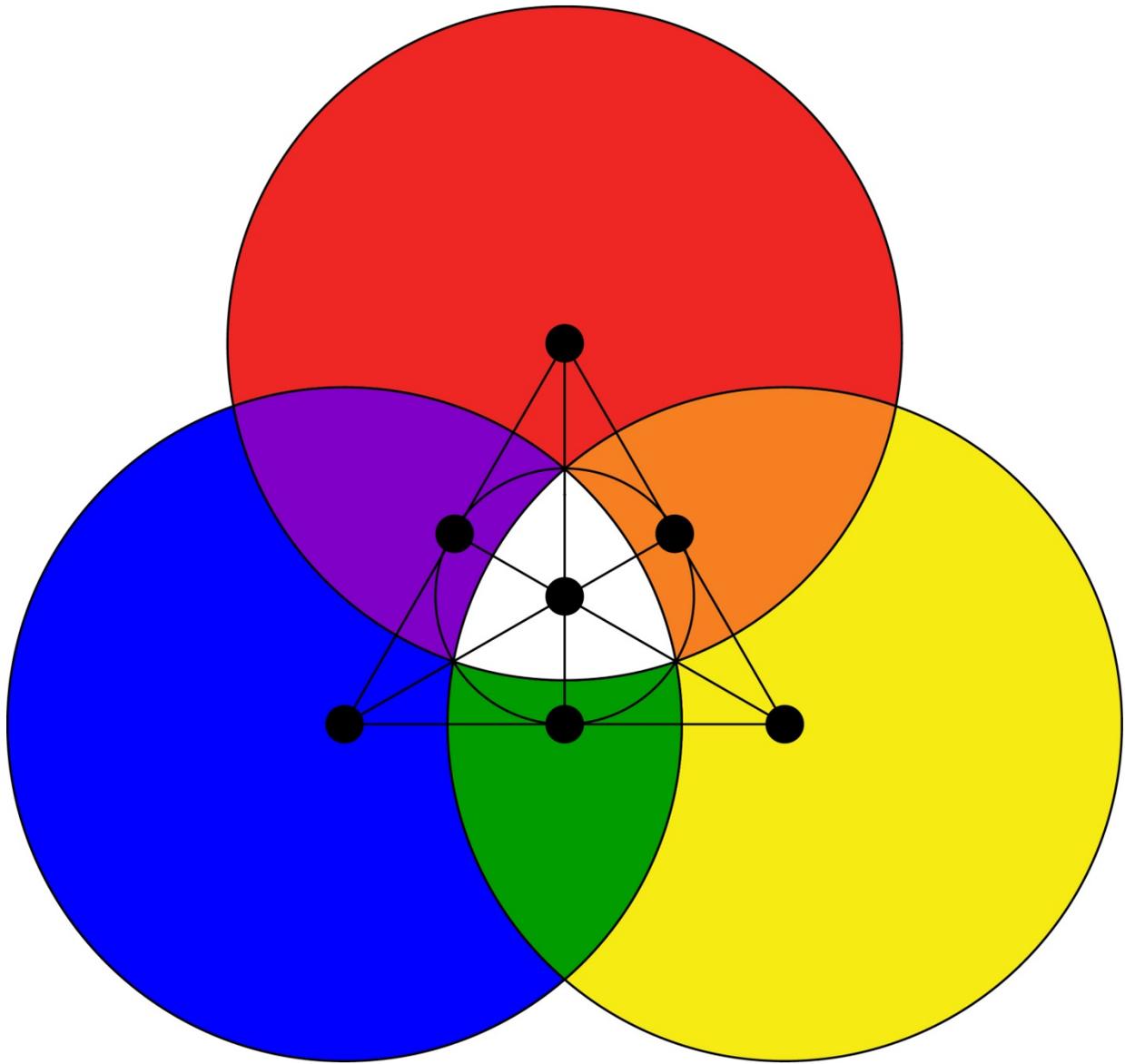

# AL-JABAR

*A Mathematical Game of Strategy*

ROBERT P. SCHNEIDER AND CYRUS HETTLE

# Concepts

The game of *Al-Jabar* is based on concepts of color-mixing familiar to most of us from childhood, and on ideas from abstract algebra, a branch of higher mathematics. Once you are familiar with the rules of the game, your intuitive notions of color lead to interesting and often counter-intuitive color combinations created in gameplay.

Because *Al-Jabar* requires some preliminary understanding of the color-mixing mechanic before playing the game, these rules are organized somewhat differently than most rulebooks. This first section details the arithmetic of combining colors used in the game. While the mathematics involved uses some elements of group theory, a foundational topic in abstract algebra, understanding this arithmetic is not difficult and requires no mathematical background. The second section explains the rules of play and how this arithmetic of colors is used in the game, and a third section describes the mathematics behind the game.

Gameplay consists of manipulating game pieces in the three primary colors red, blue and yellow, which we denote in writing by R, B, and Y respectively; the three secondary colors green, orange and purple, which we denote by G, O, and P; the color white, denoted by W; and either black or clear pieces (whichever are more readily available to players), denoted by K, which are considered to be "empty" as they are absent of color.

We refer to a game piece by its color, e.g. a red piece is referred to as "red," or R.

We use the symbol "+" to denote a combination, or grouping together, of colored game pieces, and call such a combination a "sum of colors." Any such grouping of colors will have a single color as its final result, or "sum." We use the symbol "=" to mean that two sets of pieces are equal, or interchangeable, according to the rules of the game; that is, the sets have the same sum. The order of a set of colors does not affect its sum; the pieces can be placed however you like.

Keep in mind as you read on that these equations just stand for clusters of pieces. Try to see pictures of colorful pieces, not black-and-white symbols, in your mind. Try to imagine a red piece when you read "R," a blue and a green piece when you read "B + G," and so on. In mathematics, symbols are usually just a black-and-white way to write something much prettier.

Here are four of the defining rules in *Al-Jabar*, from which the entire game follows:

$$P = R + B$$

indicates that purple is the sum of red and blue, i.e. a red and a blue may be exchanged for a purple during gameplay, and vice versa;

$$O = R + Y$$

indicates that orange is the sum of red and yellow;

$$G = B + Y$$

indicates that green is the sum of blue and yellow; and a less obvious rule

$$W = R + B + Y$$

indicates that white is the sum of red, blue and yellow, which reminds us of the fact that white light contains all the colors of the spectrum—in fact, we see in the above equation that the three secondary colors R + B, R + Y and B + Y are also contained in the sum W.

In addition, there are two rules related to the black/clear pieces. Here we use red as an example color, but the same rules apply to every color, including the black/clear piece itself:
$$R + K = R$$

indicates that a sum of colors is not changed by adding or removing a black/clear piece; and a special rule
$$R + R = K$$

indicates that two pieces of the same color (referred to as a "double") are interchangeable with a black/clear piece in gameplay. It follows from the above two rules that if we have a sum containing a double, like R + B + B, then

$$R + B + B = R + K$$

as the two blues are equal to a black/clear piece. But R + K = R so we find that

$$R + B + B = R,$$

which indicates that a sum of colors is not changed by adding or removing a double—the doubles are effectively "canceled" from the sum. It also follows from these rules that if we replace R and B with K in the above equations,

$$K + K = K$$
$$K + K + K = K, \text{ etc.}$$

We note that all groups of pieces having the same sum are interchangeable in *Al-Jabar*. For instance,
$$Y + O = Y + R + Y = R + Y + Y$$
$$= R + K = R,$$

as orange may be replaced by
$$R + Y$$

and then the double Y + Y may be canceled from the sum. But it is also true that

$$B + P = B + R + B = R + B + B$$
$$= R + K = R,$$
and even

$$G + W = B + Y + W = B + Y + R + B + Y$$
$$= R + B + B + Y + Y = R + K + K = R,$$

which uses the same rules, but takes an extra step as both G and W are replaced by primary colors.

All of these different combinations have a sum of R, so they are equal to each other, and interchangeable in gameplay:

$$Y + O = B + P = G + W = R.$$

In fact, every color in the game can be represented as the sum of two other colors in many different ways, and all these combinations adding up to the same color are interchangeable.

Every color can also be represented in many different ways as the sum of three other colors; for example
$$Y + P + G = Y + R + B + B + Y$$
$$= R + K + K = R$$
and

$$O + P + W = R + Y + R + B + R + B + Y$$
$$= R + K + K + K = R$$

are interchangeable with all of the above combinations having sum R.

An easy technique for working out the sum of a set of colors is this:

**1.** Cancel the doubles from the set;

**2.** Replace each secondary color, or white, with the sum of the appropriate primary colors;

**3.** Cancel the doubles from this larger set of colors;

**4.** Replace the remaining colors with a single piece, if possible, or repeat these steps until only one piece remains (possibly a black/clear piece). The color of this piece is the sum of the original set, as each step simplifies the set but does not affect its sum.

As you become familiar with these rules and concepts, it is often possible to skip multiple steps in your mind, and you will begin to see many possibilities for different combinations at once.

Before playing, you should be familiar with these important combinations, and prove for yourself that they are true by the rules of the game:

$$R + O = Y, \quad Y + O = R, \quad B + P = R, \quad R + P = B, \quad B + G = Y, \quad Y + G = B.$$

These show that a secondary color plus one of the primary colors composing it equals the other primary color composing it.

You should know, and prove for yourself, that

$$G + O = P, \quad O + P = G, \quad P + G = O,$$

i.e. that the sum of two secondary colors is equal to the other secondary color. You should know, and prove for yourself, that adding any two equal or interchangeable sets equals black/clear; for example
$$R + B = P \text{ so } R + B + P = K.$$

You should experiment with sums involving white—it is the most versatile color in gameplay, as it contains all of the other colors.

Play around with the colors. See what happens if you add two or three colors together; see what combinations are equal to K; take a handful of pieces at random and find its sum. Soon you will discover your own combinations, and develop your own tricks.

Included below is a table of all possible combinations of two game pieces. Combinations involving three or more pieces may be found from the table by first reducing two of the pieces to a single color. This color arithmetic is also encoded in the *Al-Jabar* logo on the cover page of this document—check it out and see what we mean.

| Color | K | R | B | Y | P | O | G | W |
|---|---|---|---|---|---|---|---|---|
| K | K | R | B | Y | P | O | G | W |
| R | R | K | P | O | B | Y | W | G |
| B | B | P | K | G | R | W | Y | O |
| Y | Y | O | G | K | W | R | B | P |
| P | P | B | R | W | K | G | O | Y |
| O | O | Y | W | R | G | K | P | B |
| G | G | W | Y | B | O | P | K | R |
| W | W | G | O | P | Y | B | R | K |

# Rules of play

**1**. *Al-Jabar* is played by 2 to 4 people. The object of the game is to finish with the fewest game pieces in one's hand, as detailed below.

**2**. One player is the dealer. The dealer draws from a bag of 70 game pieces (10 each of the colors white, red, yellow, blue, orange, green, and purple), and places 30 black/clear pieces in a location accessible to all players.

> *NOTE: Later in the game, it may happen that the black/clear pieces run out due to rule 6. In this event, players may remove black/clear pieces from the center and place them in the general supply, taking care to leave a few in the center. If there are still an insufficient number, substitutes may be used, as the number of black/clear pieces provided is not intended to be a limit.*

**3**. Each player is dealt 13 game pieces, drawn at random from the bag, which remain visible to all throughout the game.

**4**. To initiate gameplay, one colored game piece, drawn at random from the bag, and one black/clear piece are placed on the central game surface (called the "Center") by the dealer.

**5**. Beginning with the player to the left of dealer and proceeding clockwise, each player takes a turn by exchanging any combination of 1, 2 or 3 pieces from his or her hand for a set of 1, 2 or 3 pieces from the Center having an equal sum of colors.

The exception to this rule is the combination of 4 pieces R + B + Y + W, which may be exchanged for a black/clear piece. This action is called the "Spectrum" move.

> *NOTE: The shortest that a game may last is 4 moves, for a player may reduce their hand by at most 3 pieces in a turn.*

If a player having more than 3 game pieces in hand cannot make a valid move in a given turn, then he or she must draw additional pieces at random from the bag into his or her hand until a move can be made.

**6**. If a player's turn results in one or more pairs of like colors (such a pair is called a "double") occurring in the Center, then each such double is removed from the Center and discarded (or "canceled"), to be replaced by a black/clear piece.

In addition, every other player must draw the same number of black/clear pieces as are produced by cancellations in this turn.

There are two exceptions to this rule:

(*i*) Pairs of black/clear pieces are never canceled from the Center;

(*ii*) If a player's turn includes a double in the set of pieces placed from his or her hand to the Center, then the other players are not required to take black/clear pieces due to cancellations of that color, although black/clear pieces may still be drawn from cancellations of other colored pairs.

> *NOTE: The goals of a player, during his or her turn, are to exchange the largest possible number of pieces from his or her hand for the smallest number of pieces from the Center; and to create as many cancellations in the Center as possible, so as to require the other players to draw black/clear pieces.*

**7**. A player may draw additional pieces as desired at random from the bag during his or her turn.

> *NOTE: If a player finds that his or her hand is composed mostly of a few colors, or requires a certain color for a particularly effective future move, then such a draw may be a wise idea.*

**8.** A round of gameplay is complete when every player, starting with the first player, has taken a turn. Either or both of two events may signal that the game is in its final round:

(*i*) One player announces, immediately after his or her turn, that he or she has reduced his or her hand to one piece;

(*ii*) Any player having 3 or fewer pieces in hand chooses to announce, before his or her turn, that the current round will be the final round. (Otherwise, he or she may have to draw until a move can be made.)

In either case, players who have not yet taken turns in the current round are allowed to make their final moves. Even if the player who signaled the end of the game receives additional pieces in the final round, due to cancellations produced by other players, the game comes to a finish.

Drawing additional pieces is optional in the final round, and a player who cannot make a move may choose to pass without drawing.

When this final round is complete, the player with the fewest remaining pieces in hand is the winner. If two or more players are tied for the fewest number of pieces in hand, they share the victory.

A short video of gameplay may be found at [3]. Please note that in the video, clear game pieces are used to represent black/clear pieces.

## Mathematical structure

Here we outline the algebraic properties of *Al-Jabar*. This section is in no way essential for gameplay. Rather, the following notes are included to aid in understanding the mathematics behind the game and extending the game rules, which were derived using general formulas, to include sets having any number of "primary" elements, or comprised of game pieces other than colors. We encourage the interested reader to refer to [1] for more about abstract algebra, and [2] for more about mathematical games and puzzles.

The arithmetic of *Al-Jabar* in the group of the eight colors of the game is isomorphic to the addition of ordered triples in $\mathbb{Z}_2 \times \mathbb{Z}_2 \times \mathbb{Z}_2$, that is, 3-vectors whose elements lie in the congruence classes modulo 2.

The relationship becomes clear if we identify the three primary colors red, yellow, and blue with the ordered triples

$$R = (1,0,0), \ Y = (0,1,0), \ B = (0,0,1)$$

and define the black/clear color to be the identity vector

$$K = (0,0,0).$$

We identify the other colors in the game with the following ordered triples using component-wise vector addition:

$$\begin{aligned} O &= R + Y = (1,0,0) + (0,1,0) = (1,1,0) \\ G &= Y + B = (0,1,0) + (0,0,1) = (0,1,1) \\ P &= R + B = (1,0,0) + (0,0,1) = (1,0,1) \\ W &= R + Y + B = (1,0,0) + (0,1,0) + (0,0,1) = (1,1,1). \end{aligned}$$

The color-addition properties of the game follow immediately from these identities if we sum the vector entries using addition modulo 2. Then the set of colors {R, Y, B, O, G, P, W, K} is a group under the given operation of addition, for it is closed, associative, has an identity element (K), and each element has an inverse (itself).

Certain rules of gameplay were derived from general formulas, the rationale for which involved a mixture of probabilistic and strategic considerations. Using these formulas, the rules of

*Al-Jabar* can be generalized to encompass different finite cyclic groups and different numbers of primary elements, using $n$-vectors with entries in $\mathbb{Z}_m$, i.e. elements $\mathbb{Z}_m \times \mathbb{Z}_m \times \mathbb{Z}_m \times \ldots \times \mathbb{Z}_m$ ($n$ times).

In such a more general setting, there are $m$ "primary" $n$-vectors of the forms $(1,0,0,\ldots,0), (0,1,0,\ldots,0), \ldots, (0,0,\ldots,0,1)$, and the other nonzero $m$-vectors comprising the group are generated using component-wise addition modulo $m$, as above. Also, the analog to the black/clear game piece is the zero-vector $(0,0,0,\ldots,0)$.

In addition, the following numbered rules from the Rules of Play would be generalized as described here:

**2.** The initial pool of game pieces used to deal from will be composed of at least $Am^n - A$ pieces, where $A$ is at least as great as $m$ multiplied by the number of players. This pool of pieces will be divided into an equal number $A$ of every game piece color except for the black/clear or identity-element $(0,0,0,\ldots,0)$. Players will recall that the number of black/clear pieces is arbitrary and intended to be unlimited during gameplay, so this number will not be affected by the choice of $m$ and $n$.

**3.** The number of pieces initially dealt to each player will be $m^{n+1} - m - 1$.

**5.** On each turn, a player will exchange up to $n$ pieces from his or her hand for up to $n$ marbles from the Center with the same sum. The exception to this is the Spectrum, which will consist of the $n$ primary colors

$$(1,0,0,\ldots,0), (0,1,0,\ldots,0), \ldots, (0,0,\ldots,0,1)$$

together with the $n$-vector

$$(m-1, m-1, m-1, \ldots, m-1),$$

which is the generalized analog to the white game piece used in the regular game. It will be seen that these $n+1$ marbles have a sum of $(0,0,0,\ldots,0)$ or black/clear. A player must draw additional marbles if he or she has more than $n$ pieces in hand and cannot make a move.

**6.** The cancellation rule will apply to $m$-tuples (instead of doubles) of identical non-black/clear colors.

**8.** The first player to have only one piece remaining after his or her turn will signal the final round, or any player having $n$ or fewer pieces in hand may choose to do so.

Thus for the group $\mathbb{Z}_2 \times \mathbb{Z}_2 \times \mathbb{Z}_2 \times \mathbb{Z}_2$ we have $m = 2, n = 4$ and let $A = 10$. Then each player starts with 29 game pieces dealt from a bag of 10 each of the 15 non-black/clear colors, may exchange up to 4 pieces on any turn or 5 pieces in the case of a Spectrum move, and may signal the end of the game with 4 or fewer pieces in hand.

Here the Spectrum consists of the colors

$$(1,0,0,0), (0,1,0,0), (0,0,1,0), (0,0,0,1), (1,1,1,1)$$

and the cancellation rule still applies to doubles in this example, as $m = 2$.

Other cyclic groups may also be seen as sets of colors under our addition, such as

$$\mathbb{Z}_3 \times \mathbb{Z}_3 = \{(0,0), (0,1), (0,2), (1,0), (1,1), (1,2), (2,0), (2,1), (2,2)\}$$

in which every game piece either contains, for example, no red (0), light red (1) or dark red (2) in the first vector entry, and either contains no blue (0), light blue (1) or dark blue (2) in the second entry. Then we might respectively classify the nine elements above as the set

{black/clear, light blue, dark blue, light red, light purple, bluish purple,
dark red, reddish purple, dark purple}.

Of course, other colors rather than shades of red and blue may be used, or even appropriately selected non-colored game pieces.

Further generalizations of the game rules may be possible—for instance, using $n$-vectors in $\mathbb{Z}_{a_1} \times \mathbb{Z}_{a_2} \times \ldots \times \mathbb{Z}_{a_n}$ where the subscripts $a_i$ are not all equal—and new games might be produced by other modifications to the rules of play or the game pieces used.

As noted above, the *Al-Jabar* logo encodes the algebra of the game. Each node on the Fano plane diagram represents the color on which it falls; the sum of any two nodes lying on the same line segment (or on the inner circle) is equal to the third node on that segment, and all three points on a given segment (or the inner circle) taken together are equal to black/clear.

## Acknowledgments

The authors are grateful to the following people and organizations: Maxwell Schneider, who inspired the game's conception; Prof. Paul Eakin (University of Kentucky) whose Modern Algebra course motivated the structure of the game; artist Robert Beatty for rendering the game logo; David Unzicker, Barbara Hettle, Marci Schneider, Trish Smith, Prof. Ron Taylor and his students (Berry College), Prof. Colm Mulcahy and his students (Spelman College), Prof. Neil Calkin (Clemson University), Adam Jobson (University of Louisville), and the Lexington Board Game Group (Lexington, Kentucky) for play-testing, support and suggestions; and to Gathering for Gardner who first published the game rules in honor of Martin Gardner's birthday (Celebration of Mind, October 21, 2011).